\documentclass{article}    

\usepackage{amsmath}
\usepackage{amsfonts}
\usepackage{amsthm}
\usepackage{amssymb}
\usepackage{graphicx}

\title{\bf Counting Mountain-Valley Assignments for Flat Folds\footnote{This paper appeared in {\em Ars Combinatoria}, Vol. 67 (2003), 175--188.}}
\author{Thomas Hull\\
Department of Mathematics\\
Merrimack College\\
North Andover, MA 01845}

\date{}

\newtheorem{thm}{Theorem}

\newtheorem{lemma}{Lemma}

\begin{document} 

\maketitle

\begin{abstract}
We develop a combinatorial model of paperfolding for the purposes of
enumeration.  A planar embedding of a graph is called a {\em crease pattern}
if it represents the crease lines needed to fold a piece of paper into something.
A {\em flat fold} is a crease pattern which lies flat when folded, i.e. can be
pressed in a book without crumpling.  Given a crease pattern $C=(V,E)$, a {\em
mountain-valley (MV) assignment} is a function $f:E\rightarrow \{$M,V$\}$ which
indicates which crease lines are convex and which are concave, respectively. A
MV assignment is {\em valid} if it doesn't force the paper to self-intersect
when folded.  We examine the problem of counting the number of valid MV assignments
for a given crease pattern.  In particular we develop recursive functions that count
the number of valid MV assignments for {\em flat vertex folds}, crease patterns with
only one vertex in the interior of the paper.  We also provide examples,
especially those of Justin, that  illustrate the difficulty of the general
multivertex case.
\end{abstract}

\section{Introduction}

The study of origami, the art and process of paperfolding, includes many
interesting geometric and combinatorial problems.  (See \cite{Hull} and
\cite{Justin} for more background.)  In origami mathematics, a {\em fold}
refers to any folded paper object, independent of the number of folds done
in sequence.  The {\em crease
pattern} of a fold is a planar embedding of a graph which represents the creases that
are used in the final folded
object.  (This can be thought of as a structural blueprint of the fold.)
Creases come in two types: {\em mountain
creases}, which are convex, and {\em valley creases}, which are concave.
Clearly the type of a crease depends on which side of the paper we look at,
and so we assume we are always looking at the same side of the paper.
In this paper
we will concern ourselves with the following question about {\em flat
folds}, i.e., origami that can, when completed, be pressed in a book without
crumpling: 

\begin{quote}
Given a crease pattern that can fold flat, how many different ways can we
assign
mountain and valley creases and still collapse it?
\end{quote}

More formally, we define a {\bf MV assignment} of a given a crease pattern $C=(V,E)$
to be a function
$f:E\rightarrow\{M,V\}$.  MV assignments that can actually be folded are called {\em
valid}, while those which do not admit a flat folding (i.e. force the paper to
self-intersect in some way) are called {\em invalid}. A complete answer to the
problem of counting the number of valid MV assignments of a given crease pattern is
currently inaccessible.  Any given crease pattern can be collapsed in many
different ways.  The purpose of this paper it to present, formalize, and expand the
known results for counting the number of valid MV assignments for a given flat fold
crease pattern, focusing primarily on the single vertex case.  At the same time we
will discover that many of the results that hold for folding a sheet of paper flat
also hold for folding {\em cone shaped} paper (that has less than
$360^\circ$ around a vertex) with a crease pattern whose only vertex is at the apex
of the cone.

\section{Preliminaries}

Whether or not a crease pattern will fold flat is not completely determined by
a MV assignment; other factors come into play including the arrangement of the
layers of paper and whether or not this arrangement will force the paper to intersect
itself when folded, which is {\em not} allowed.  We present a few basic Theorems
relating to necessary and sufficient conditions for flat-foldability. These Theorems
appear in their cited references without proof.  While Kawasaki, Maekawa, and
Justin undoubtedly had proofs of their own, the proofs presented below were devised
by the author and Jan Siwanowicz at the 1993 Hampshire College Summer Studies in
Mathematics.

\begin{thm}[Kawasaki-Justin \cite{Justin1},
\cite{Justin2}, \cite{Kawa3}]\label{kj}
Let $v$ be a vertex of degree $2n$ in an origami crease pattern of
and let
$\alpha_1, ..., \alpha_{2n}$ be the consecutive angles between the creases.  Then the
creases adjacent to $v$ will (locally) fold flat if and only if
\begin{equation}\label{iso}
\alpha_1-\alpha_2+\alpha_3-\cdots -\alpha_{2n}=0.
\end{equation}
\end{thm}

\noindent{\bf Proof:} The equation easily follows by 
considering a simple closed curve
which winds around the vertex.  This curve mimics the path of an ant
walking around the vertex on the surface of the paper after it is folded.  We
measure the
angles the ant crosses as positive in one direction and negative in the other. 
Arriving at the point where the ant started means that this alternating sum is
zero. The converse is left as an exercise. (See \cite{Hull}.) $\qed$

\begin{thm}[Maekawa-Justin \cite{Justin2}, \cite{conn}]\label{mj}
Let $M$ be the number of mountain creases and $V$ be the number of valley
creases adjacent to a vertex in a flat origami crease pattern.  Then $M-V=\pm 2$.
\end{thm}

\noindent{\bf Proof:} (Jan Siwanowicz) If $n$ is the number of creases, then $n=M+V$. 
Fold the paper flat and 
consider the cross-section
obtained by clipping the area near the vertex from the paper; the cross-section
forms a flat polygon.  If we view each interior
$0^\circ$ angle as a valley crease and each interior $360^\circ$ angle as a mountain
crease, then $0V+360M=(n-2)180=(M+V-2)180$, which gives $M-V=-2$. On the other hand,
if we view each
$0^\circ$ angle as a \emph{mountain} crease and each $360^\circ$ angle as a
\emph{valley} crease (this corresponds to flipping the paper over), then
we get $M-V=2$.  $\qed$

\vspace{.1in}

We refer to Theorems \ref{kj} and \ref{mj} as the K-J Theorem and the M-J Theorem,
respectively. Justin \cite{Justin}
refers to equation (\ref{iso}) as the {\em isometries condition}. 
The K-J Theorem is sometimes
stated in the equivalent form that the sum of every other angle around $v$ equals
$180^\circ$, but this is only true if the vertex is on a flat sheet of paper. 
Indeed, notice that the proofs of the K-J and M-J Theorems do not use the fact that
$\sum \alpha_i = 360^\circ$.  Thus these two theorems are also valid for flat vertex
folds where $v$ is at the apex of a cone-shaped piece of paper.  We will require
this generalization later.

Note that while the K-J Theorem does assume that the vertex has even degree, the M-J
Theorem does not.  Indeed, the M-J Theorem can be used to prove this fact.  Let
$v$ be a vertex in a crease pattern that folds flat and let $n$ be the degree of $v$.
Then $n=M+V=M-V+2V=\pm 2 + 2V$, which is even.

In their present form neither of these theorems generalize to handle more than one
vertex in a crease pattern.\footnote{Although Kawasaki has been able to reformulate
his Theorem to say something about flat origami crease patterns in general (see
\cite{Kaw1}) and Justin posits necessary and sufficient conditions for global
flat-foldability (see \cite{Justin}), we won't be using these results here.}
To illustrate the difficulty involved in determining the number of valid
MV assignments in a flat multiple vertex fold, we present an exercise, which the
reader is encouraged to attempt. 

\vspace{.2in}

\begin{figure}[h]
\centerline{\includegraphics[scale=.55]{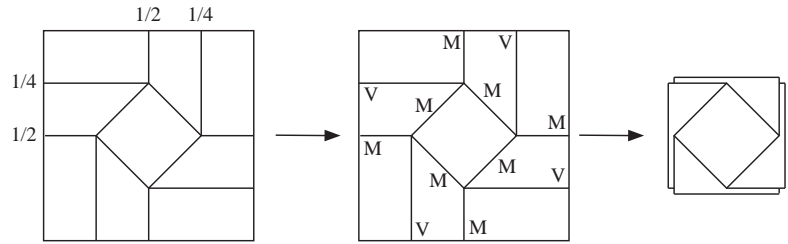}}
\caption{a square twist}\label{e1}
\end{figure}

\noindent{\bf Exercise:} Figure \ref{e1} displays the crease pattern for an origami
fold called a {\em square twist}, together with a valid MV assignment. Of the
$2^{12}$ different possible MV assignments for this crease pattern, only 16 are
valid.  Can you find them all?

\begin{figure}[h]
\centerline{\includegraphics[scale=.7]{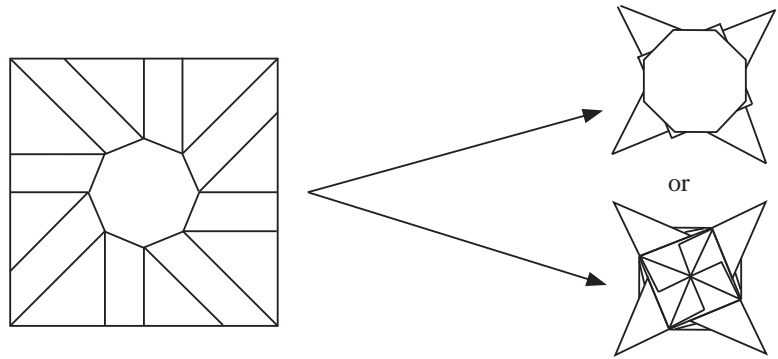}}
\caption{an octagon twist}\label{e2}
\end{figure}

In contrast, Figure \ref{e2} shows an {\em octagon twist}, which is
a bit more difficult to fold than the square twist.  Although this has more creases
than the previous exercise, this octagon twist has only {\em two}
different valid MV assignments.  Indeed, experimentation with this crease pattern
makes it apparent that the inner octagon must be all mountain creases or all valley
creases, which forces the assignment of the remaining creases. (If, however, the
octagon is made to be larger relative to the paper's boundary, then more valid MV
assignments can be possible.)

Note that when considering crease patterns for origami folds, we {\em do not} include the boundary
of the paper in the crease pattern, since no folding is actually taking place there.  Thus we {\em
ignore} vertices of the crease pattern on the boundary of the paper, since any results (like the K-J
and M-K Theorems) will not apply to these vertices.  Only vertices in the interior of the paper are
considered, which is natural for the present study because we'll be primarily investigating local
properties.  (I.e., how the paper behaves around a single vertex.)

\section{Flat vertex folds}
We now restrict ourselves to {\em flat
vertex folds}, which are folds whose crease patterns contain only one vertex in the
interior of the paper.  
For our purposes we may consider a vertex $v$ to be completely determined by the angles
between its crease lines, and thus we write $v=(\alpha_1,...,\alpha_{2n})$ where the
$\alpha_i$ denote the consecutive angles around $v$.
Throughout 
we let
$l_1,..., l_{2n}$ denote the creases adjacent to a vertex
$v$ where $\alpha_i$
is the angle between creases $l_i$ and $l_{i+1}$ for $i=1,...,2n-1$ and $\alpha_{2n}$
is between $l_{2n}$ and $l_1$.  

Let us consider a flat vertex fold, which we know must satisfy the M-J Theorem. Given
a specific MV assignment to the creases, if $M-V=2$ then we say that
$v$ {\em points up}. If $M-V=-2$ then  we say it
{\em points down}. Because a vertex that points up can be made to point down by
reversing all the mountain and valley creases, and vice-versa, we will study all the
possible valid MV assignments of the crease pattern by considering only
cases where
$v$ points up, knowing that each of these cases has a pointing down counterpart.

In general, we denote
\begin{eqnarray*}
C(\alpha_1,...,\alpha_{2n}) & = & \mbox{the number of valid MV
assignments of a flat}\\
& & \mbox{vertex fold with consecutive angles $\alpha_1,...,\alpha_{2n}$.}
\end{eqnarray*}

Note that by the above observation $C(\alpha_1,...,\alpha_{2n})$ 
is always even.

\begin{figure}
\centerline{\includegraphics[scale=.9]{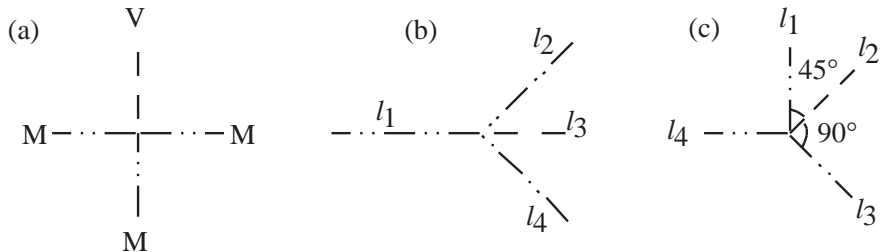}}
\caption{(a) $C(\alpha_1,...,\alpha_4)=8$, (b) $C(\alpha_1,...,\alpha_4)=6$, (c)
$C(\alpha_1,...,\alpha_4)=4$.}\label{fig2}
\end{figure}

Let us consider some basic examples for computing $C(\alpha_1,...,\alpha_{2n})$. Suppose 
the vertex has 4
creases.  We will demonstrate that
$C(\alpha_1,...,\alpha_4)$ can take on the values 8, 6, or 4 in this case, depending on the angles
between the creases.  Examine the three flat vertex folds with MV assignments shown
in Figure \ref{fig2}. (We follow standard origami notation by denoting mountain and
valley creases by two different kinds of dashed lines.)

To see how the examples in Figure \ref{fig2} work, assume that the vertices are pointing up. 
Thus in all three cases we'll have 3 mountains and 1 valley.  In (a) the single
valley crease could be any of the four crease lines since the angles are all the
same.  (Just fold a square piece of paper in half, then in half again.) This gives
${4\choose 1}=4$ possibilities, and each of these has a pointing down counterpart. 
Thus
$C(\alpha_1,...,\alpha_4)=8$.

In (b) notice that crease line $l_1$ cannot be the valley.  (If it were, then when
we try to fold it flat the regions of paper bordering $l_1$ would have to intersect
crease $l_2$ because the acute angles around $l_2$ can't completely contain the
obtuse angles around $l_1$. The reader is encouraged to experiment.)  Thus
there are only 3 positions where the valley crease can be assigned, giving
$C(\alpha_1,...,\alpha_4)=6$.

In (c), notice that either $l_1$ or $l_2$ must be the valley.  This is because $l_1$
and $l_2$ make a $45^\circ$ angle, with $90^\circ$ angles on the left and right of
it.  In this ``big angle-little angle-big angle" case, creases
$l_1$ and $l_2$ can't both have the same MV parity because that would force the two
$90^\circ$ angles to cover up the smaller $45^\circ$ angle on the same side of the
paper, causing the regions of paper made by the two large angles to intersect one
another.  Thus in this crease configuration the valley crease has only  the
possibilities
$l_1$ or
$l_2$, giving
$C(\alpha_1,...,\alpha_4)=4$.

Notice that the situation in example (c) can be generalized.  In particular, when we
have an angle in a flat vertex fold, we might get into trouble if we make the
creases bordering the angle be both mountains or both valleys, but making the
creases different will always work.  We will need to use this later, so we state it
explicitly.

\vspace{.1in}

\noindent{\bf Observation 1:} In a flat vertex fold, if we have consecutive
angles $\alpha_{i-1}, \alpha_i$ and $\alpha_{i+1}$, then we can
always assign $l_i$ to be a mountain and $l_{i+1}$ to be a valley, or vice-versa,
without risk of a forced self-intersection of the paper among the parts of the paper
made by these angles.

\vspace{.1in}

Notice further that example (a) in Figure \ref{fig2}, where all the angles between the
creases were equal, gave us the most variability.  This is true in general; if all
the angles between the creases are equal then, by symmetry, it doesn't matter where
the valley creases are placed.  Thus if we have $2n$ creases and the vertex points
up, then any
$n-1$ of them can be valleys to satisfy the M-J Theorem.  The number of ways we can
choose these valley creases is ${2n\choose n-1}$, and each of these has a pointing
down counterpart. This gives us the upper bound
$C(\alpha_1,...,\alpha_{2n})\leq 2{2n\choose n-1}$.  Furthermore, this would also hold if our flat vertex
fold was at the apex of a cone-shaped piece of paper.  We have proven half of the
following:

\begin{thm}\label{ulb}
Let $v=(\alpha_1,...,\alpha_{2n})$ be the vertex in a flat vertex fold, on either a
flat piece of paper or a cone.  Then
$$2^n\leq C(\alpha_1,...,\alpha_{2n})\leq 2{2n\choose n-1}$$
are sharp bounds.
\end{thm}

A number of people have discovered the lower bound in Theorem \ref{ulb}.  Azuma
\cite{Azuma} presented the result without proof, Justin \cite{Justin} provides all
the elements of a proof but does not state the result explicitly, and Ewins and Hull
\cite{EH} independently constructed the proof given below.

\vspace{.1in}

\noindent{\bf Proof of the lower bound:} 
Imagine we have a flat vertex fold on a flat piece of paper or a cone and suppose
$\alpha_i$ is the smallest angle surrounding the vertex $v$.  (Or one of the
smallest, if there is a tie.)  If $l_i$ and $l_{i+1}$ are the creases on the left
and right of angle
$\alpha_i$, then by Observation 1 we have at {\em least two} possibilities
for the MV assignment of $l_i$ and $l_{i+1}$. $f(l_i, l_{i+1})$ could
be (M,V) or (V,M).  (Of course, there might be other possibilities.)  

Thus if we fold $l_i$ and $l_{i+1}$ using one of these two possibilities and 
{\em fuse}, or identify the layers of paper together, then the paper will turn into a
cone (unless it already is a cone, in which case it will become a smaller cone) and
angles
$\alpha_{i-1},
\alpha_i,
\alpha_{i+1}$ will become a new angle with measure
$\alpha_{i-1}-\alpha_i+\alpha_{i+1}$, which will be positive because $\alpha_i$ was
one of the smallest angles.  Since the original flat vertex fold can fold flat, our
new cone will also fold flat along the remaining crease lines
$l_1, l_2,... l_{i-1}, l_{i+2},...,l_{2n}$.  

In other words, we can repeat this process. Take the smallest angle in our new cone,
fold its bordering creases in one of the two guaranteed possible ways ((M,V) or
(V,M)), fuse them together, and repeat.  Each time we do this we eliminate two
creases and count at least two possible mountain-valley configurations for those
creases.  Eventually there will only be two creases left in our cone, and these
can either be both mountains or both valleys, by the M-J
Theorem.  If we started with
$2n$ creases, we'll have eliminated a total of $n$ pairs of creases, with at least
two MV assignment choices per pair, giving us
$C(\alpha_1,...,\alpha_{2n})\geq 2^n$. $\qed$

This lower bound becomes equality ($C(\alpha_1,...,\alpha_{2n})=2^n$) for {\em generic} flat vertex
folds, which are those in which the angles are chosen so that none are
consecutively equal and none of the combined angles are equal to their neighbors
throughout the recursive process outlined above.  For example, if we have six
creases with angles
$100^\circ, 70^\circ, 50^\circ, 40^\circ, 30^\circ, 70^\circ$ surrounding a vertex
$v$, then we have $C(\alpha_1,...,\alpha_{2n})=2^3=8$.

In any case, we see that a simple formula for $C(\alpha_1,...,\alpha_{2n})$ in terms of $n$ alone is not
possible.  To actually compute $C(\alpha_1,...,\alpha_{2n})$ more information, in particular the
values of the angles between the creases, is needed.

\subsection{Many equal angles in a row}

It will be useful for us to introduce the following notation:  If $l_i,...,l_{i+k}$
are consecutive crease lines in a flat vertex fold which have been given a MV
assignment, let $M_{i,...,i+k}=$ the number of mountains and $V_{i,...,i+k}=$ the
number of valleys among these crease lines.

Suppose that somewhere in our flat vertex fold (in either a flat piece
of paper or a cone) we have
$\alpha_i=\alpha_{i+1}=\alpha_{i+2}=\cdots =\alpha_{i+k}$ and $\alpha_{i-1}>
\alpha_i$ and
$\alpha_{i+k+1}>\alpha_{i+k}$.  (Note that if $\alpha_1$ and $\alpha_{2n}$ appear
in our sequence of equal angles, we may relabel so that they do not.)
\footnote{Also note that the case where $k=2n-2$ is impossible.  Indeed, this would imply
that $i=1$ and we have $\alpha_1=\cdots=\alpha_{2n-1}$ and $\alpha_{2n}$ is bigger than all the other
angles.   But then the K-J Theorem implies that $n\alpha_1=(n-1)\alpha_1+\alpha_{2n}$, or
$\alpha_1=\alpha_{2n}$, a contradiction.} If
$k=0$ then we have a large angle, then a small one, then a large one, which is the
same situation that we saw in the example in Figure \ref{fig2} (c), above.  Thus we get that creases
$l_i$ and
$l_{i+1}$ cannot both be valleys or both be mountains.  That is, $M_{i,i+1}-
V_{i,i+1}=0$. If $k>0$, then we have several consecutive angles of the same
measure, and there will be many more possibilities for MV assignments.
The following Theorem presents the general result.  (Note that while \cite{Justin}, 
\cite{koehler}, and \cite{lunnon1} do not state this result explicitly, Justin's
work on flat foldings, and Lunnon and Koehler's work on folding and arranging
postage stamp arrays is similar enough to make it clear that this result was known to
those authors.)

\begin{thm}\label{angles1}
Let $v=(\alpha_1,...,\alpha_{2n})$ be a flat vertex fold in either a piece of paper or a cone, and
suppose we have
$\alpha_i= \alpha_{i+1}=\alpha_{i+2}=\cdots =\alpha_{i+k}$ and $\alpha_{i-1}>
\alpha_i$ and
$\alpha_{i+k+1}>\alpha_{i+k}$ for some
$i$ and $k$.  Then
$$M_{i,...,i+k+1}-V_{i,...,i+k+1} = \left\{\begin{array}{cl}
0 & \mbox{if $k$ is even}\\
\pm 1 & \mbox{if $k$ is odd.}\end{array}\right.$$
\end{thm}

\begin{figure}
\centerline{\includegraphics[scale=.8]{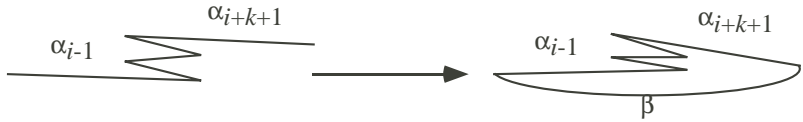}}
\caption{$k$ is even.}\label{fig3}
\end{figure}

\noindent{\bf Proof:} The result follows from an application of the M-J
Theorem.  If $k$ is even, then the cross-section of the paper around the
creases in question might look as shown in Figure \ref{fig3}.\footnote{We say ``might" because
the equal angles may be twisted among themselves in a number of different ways.} If
we consider this sequence of angles by itself and imagine adding a section of paper
with angle $\beta$ to connect the loose ends at the left and right (see
Figure \ref{fig3}), then we'll have a flat-folded cone which must satisfy the  
the M-J Theorem.  The angle $\beta$ added two extra creases, both of which
must be mountains or both valleys.  We may assume that the vertex points up, and thus
we subtract two from the result of the M-J Theorem to get
$M_{i,...,i+k+1}-V_{i,...,i+k+1}=0$.

\begin{figure}
\centerline{\includegraphics[scale=.65]{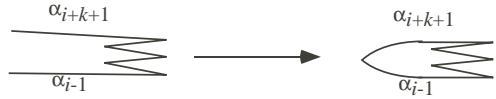}}
\caption{$k$ is odd.}\label{fig4}
\end{figure}

If $k$ is odd (Figure \ref{fig4}), then this angle sequence, if considered by itself, will
have the loose ends from angles $\alpha_{i-1}$ and $\alpha_{i+k+1}$ pointing in the
same direction.  If we glue these together, possibly extending one of them if
$\alpha_{i-1}\not=\alpha_{i+k+1}$, then the M-J Theorem may be applied.  After
subtracting (or adding) one to the result of the M-J Theorem because of the extra
crease made when gluing the loose flaps, we get $M_{i,...,i+k+1}-V_{i,...,i+k+1}=\pm
1$.  $\qed$

\vspace{.1in}

In \cite{Justin}, Justin uses this result to convert a flat vertex fold
into a {\em circular word} with parentheses to denote where Theorem \ref{angles1}
can be iteratively applied.  This provides a mechanism to enumerate all the valid
MV assignments for the creases.  A similar strategy would be to create
recursive formulas from Theorem \ref{angles1}. 

\begin{thm}\label{recur}
Let $v=(\alpha_1,...,\alpha_{2n})$ be a flat vertex fold in either a piece of paper or a cone, and
suppose we have
$\alpha_i= \alpha_{i+1}=\alpha_{i+2}=\cdots =\alpha_{i+k}$ and $\alpha_{i-1}>
\alpha_i$ and
$\alpha_{i+k+1}>\alpha_{i+k}$ for some
$i$ and $k$.  Then
$$C(\alpha_1,...,\alpha_{2n}) = 
{k+2\choose
{k+2\over 2}}C(\alpha_1,...,\alpha_{i-2},\alpha_{i-1}-\alpha_i+\alpha_{i+k+1},
\alpha_{i+k+2},..., \alpha_{2n})$$
if $k$ is even, and
$$C(\alpha_1,...,\alpha_{2n}) =
{k+2\choose {k+1\over
2}}C(\alpha_1,...,\alpha_{i-1},\alpha_{i+k+1}, ...,
\alpha_{2n})$$
if $k$ is odd.
\end{thm}

\noindent{\bf Proof:} 
If $k$ is even, then Theorem \ref{angles1} gives us
$M_{i,...,i+k+1}-V_{i,...,i+k+1}=0$, which means among the $k+2$ creases
$l_i,...,l_{i+k+1}$, any $(k+2)/2$ of them can be valleys, and the rest mountains,
since all the angles are the same.  If we take one of these possibilities and fuse
the layers of paper around these angles together, then angles
$\alpha_{i-1},..., \alpha_{i+k+1}$ will be replaced with one angle with measure
$\alpha_{i-1}-\alpha_i+\alpha_{i+k+1}$.  This gives us the stated recursion.

If $k$ is odd, then $M_{i,...,i+k+1}-V_{i,...,i+k+1}=\pm 1$.  Thus we could pick any $(k+1)/2$ of
the $k+2$ creases $l_i,...,l_{i+k+1}$ to be mountains and the rest valleys, or vice-versa.  
Thus there are $2{k+2\choose
(k+1)/2)}$ MV-assignments for these creases.
However, because
$k$ is odd, fusing all these layers together will create a new crease line whose
mountain-valley assignment will be forced and ruin our hopes of recursion.  To avoid
this, we allow one of the crease lines to remain {\em unassigned} and divide the number
of MV-assignments by two.  When the folded layers of paper are fused together, the
angles $\alpha_i, ..., \alpha_{i+k}$ will be absorbed by the angles $\alpha_{i-1}$ or
$\alpha_{i+k+1}$, which gives the stated recursion. $\qed$

\section{Examples illustrating the utility of Theorem \ref{recur}}

Theorem \ref{recur} provides us with a very efficient algorithm for computing 
\linebreak $C(\alpha_1,...,\alpha_{2n})$ 
for any flat vertex fold $v$. Examine the smallest angle.  Its neighbors will either
be larger than or equal to it, and thus we'll have an angle sequence satisfying the
conditions of Theorem \ref{recur}.  Repeat this with the new collection of angles,
until all the angles are equal.  Then the upper bound from Theorem \ref{ulb} can be
applied.

\vspace{.1in}

\noindent{\bf Example 1:} An earlier example that achieved the lower bound formula
in Theorem \ref{ulb} had six crease lines with angles $100^\circ, 70^\circ, 50^\circ,
40^\circ, 30^\circ, 70^\circ$.  Applying Theorem \ref{recur} recursively yields
\begin{eqnarray*}
C(100, 70, 50, 40, 30, 70) & = & {2\choose 1}C(100, 70, 50, 80)\\
& = & {2\choose 1}{2\choose 1}C(100, 100) \\
& = & {2\choose 1}{2\choose 1}2 = 8.
\end{eqnarray*}
{\bf Example 2:} In \cite{Justin} Justin gives the following example with eight
crease lines: $20^\circ$, $10^\circ$, $40^\circ$, $50^\circ$, $60^\circ$, $60^\circ$,
$60^\circ$, $60^\circ$.  Here we find that
\begin{eqnarray*}
C(20,10,40,50,60,60,60,60) & = & {2\choose 1}C(50, 50, 60,60,60,60)\\
& = & {2\choose 1}{3\choose 1}C(60,60,60,60)\\
& = & {2\choose 1}{3\choose 1}2{4\choose 1} =48.
\end{eqnarray*}

Notice that when we have eight creases, Theorem \ref{ulb} only tells us that $C(\alpha_1,...,\alpha_{2n})$
is anywhere between 16 and 112.  Theorem \ref{recur}, however, gives the exact value
of $C(\alpha_1,...,\alpha_{2n})$ with at most 4 computations.

\section{Multiple vertex folds}

Counting valid MV assignments of crease patterns with more than one vertex can be
very difficult.  To illustrate this, we examine a deceptively simple class of
origami folds with more than one vertex: flat origami folds whose
crease patterns are just equally-spaced grids of perpendicular lines.  This would be
a ``fold" where the paper gets folded up into a small square.

Note that  Koehler
\cite{koehler} and Lunnon \cite{lunnon1}, \cite{lunnon2}, among others, tackled what
is known as the {\bf postage stamp problem} or {\bf map-folding problem}. Here one is
given an
$m\times n$ array (sheet) of equal-sized postage stamps and the problem is to count
the number of ways one can fold them up, independent of the MV assignment.  The
authors listed above give complicated algorithms for doing this, especially for the
case where
$n=1$ and we have a strip of stamps.  However, their approach also counts the number
of different ways one can {\em arrange the layers}.  This is not, therefore, the same
as counting the number of valid MV assignments.

\begin{figure}
\centerline{\includegraphics[scale=.7]{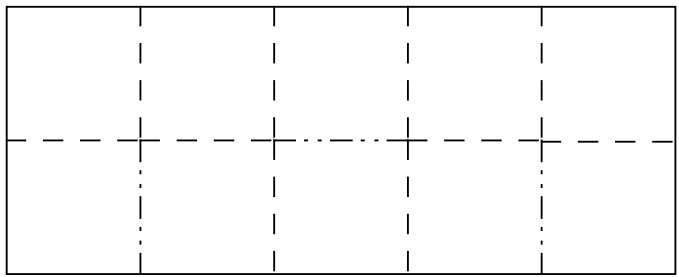}}
\caption{An impossible MV assignment for a $2\times 5$ array of stamps.}\label{stamp}
\end{figure}

Let $S_{m,n}$ be an $m\times n$ array of equal-sized postage stamps and let
$C(S_{m,n})$ denote the number of different valid MV assignments that
will fold $S_{m,n}$ into a single stamp-sized pile.  
This crease pattern
will be an $(m-1)\times (n-1)$ lattice of vertices, each of degree four.  Start with
the upper-left vertex, which has 8 different possible MV assignments. The next one to
the right then has only 4, since one of its creases is already set (that crease had
two possibilities, so this divides the number of possibilities for this vertex in
half, giving 4).  Continuing to move to the right, we have 4 possible MV assignments
for each of the remaining vertices in the top row.  The first vertex in the second
row will also have 4 possible MV assignments, but the rest of the vertices in that
row will have only 2 possible MV assignments, since their left creases and top
creases are already set.  The same will be true for the third row as well as the
remaining rows.  Thus the total number of valid MV assignments for
$S_{m,n}$ can be bounded:
$$C(S_{m,n}) \leq  8\cdot 4^{m-2}\cdot 4^{n-2}\cdot 2^{(m-2)(n-2)} =  2^{mn-1}$$
Equality is not always achieved because not all of these MV assignments are valid.
Justin in
\cite{Justin} gives a number of impossible mountain-valley assignments for the case
when
$n=2$ and $m=5,6$ and $7$.  One of the simplest is shown in Figure \ref{stamp}.

The reader is highly encouraged to try folding this nefarious crease pattern. Thus,
we can prove that
$C(S_{m,n})\leq 2^{mn-1}$, and it can be shown that we get equality when $n=1$ or $m$
and
$n$ are both less than 5.  Obtaining a better formula is an open problem.  

\section{Conclusion}

Theorem \ref{recur} provides us with a linear-time algorithm for
computing the number of valid MV assignments that can be used on a given flat vertex
fold.  The example given in Figure 5 is enough to illustrate how the equivalent
problem for flat multiple vertex folds is very daunting, indeed.  Several questions
present themselves to those who would like to work further in this area:
Might the results presented here be extended in some way to flat folds with two or
three vertices?  Are there families of flat multiple vertex folds $F$ for which
$C(F)$ is easy to compute?  

\vspace{.2in}

\noindent{\bf Acknowledgements:} The author would like to thank sarah-marie belcastro,
Michael Bradley, Jacques Justin, Jim Tanton, and an anonymous referee for valuable suggestions on
earlier drafts of this paper.

\end{document}